\documentclass[aos,preprint]{imsart}
\RequirePackage[OT1]{fontenc} \RequirePackage{amsthm,amsmath}
\RequirePackage[numbers]{natbib}
\RequirePackage[colorlinks,citecolor=blue,urlcolor=blue]{hyperref}

\usepackage{color,latexsym,amsfonts,amssymb}
\newcommand{\blue}[1]{\textcolor{blue}{#1}}

\newcommand{\red}[1]{\textcolor{red}{#1}}
\usepackage{amsfonts}
\usepackage{amssymb}

\arxiv{math.PR/0000000} \startlocaldefs
\numberwithin{equation}{section} \theoremstyle{plain}

\endlocaldefs
\newtheorem{lemma}{Lemma}[section]

\newtheorem{theorem}{Theorem}[section]

\newtheorem{example}{Example}[section]

\begin{document}

\newcommand{\tY}{\tilde{Y}}
\newcommand{\tW}{\tilde{W}}

\newcommand{\W}{W^{*}}
\newcommand{\mathds}{\mathbb}
\newcommand{\txi}{\tilde{\xi}}

\newcommand{\ul}[1]{{\bf #1}}
\newcommand{\vv}{\vspace{8mm}}
\newcommand{\vs}[1]{\vspace{#1cm}}
\newcommand{\noi}{\noindent}
\newcommand{\ds}{\displaystyle}
\newcommand{\lra}{\longrightarrow}
\newcommand{\la}{\lambda}
\newcommand{\lrla}{\Longleftrightarrow}
\newcommand{\beq}{\begin{eqnarray*}}
\newcommand{\eeq}{\end{eqnarray*}}
\newcommand{\beqn}{\begin{eqnarray}}
\newcommand{\eeqn}{\end{eqnarray}}
\newcommand{\ta}{\tau}
\newcommand{\ra}{\rightarrow}
\newcommand{\mb}{\mbox}
\newcommand{\vp}{\varepsilon}
\newcommand{\ep}{\epsilon}
\newcommand{\al}{\alpha}
\newcommand{\var}{\mbox{Var}}
\newcommand{\bi}{\begin{itemize}}
\newcommand{\ei}{\end{itemize}}
\newcommand{\sta}{\stackrel}
\newcommand{\La}{\Lambda}
\newcommand{\be}{\begin{equation}}
\newcommand{\ee}{\end{equation}}
\newcommand{\nn}{\nonumber}

\newcommand{\lbl}{\label}
\newcommand{\barB}{\bar{B}}
\newcommand{\barb}{\bar{b}}
\newcommand{\bars}{\bar{S}}
\newcommand{\barv}{\bar{V}}
\newcommand{\barx}{\bar{X}}
\newcommand{\bbox}{\nobreak\quad\vrule width4pt depth2pt height4pt}
\newcommand{\eq}[1]{$(\ref{#1})$}

\newcommand{\ignore}[1]{}{}
\newcommand{\f}{\frac}
\newcommand{\wt}{\widetilde}
\renewcommand{\(}{\left(}
\renewcommand{\)}{\right)}
\newcommand{\Xn}{X_1,\dots,X_n}
\newcommand{\Yn}{Y_1,\dots,Y_n}
\newcommand{\sn}{\sum_{i=1}^n}
\newcommand{\C}{{\bf C}}
\renewcommand{\P}{{\bf P}}
\newcommand{\Z}{{\bf Z}}
\newcommand{\R}{{\bf R}}
\newcommand{\E}{{\bf E}}
\newcommand{\de}{\delta}
\newcommand{\De}{\Delta}
\newcommand{\ga}{\gamma}
\newcommand{\Ga}{\Gamma}
\newcommand{\si}{\sigma}
\newcommand{\Si}{\Sigma}
\newcommand{\Th}{\Theta}
\newcommand{\pa}{\parallel}
\newcommand{\ka}{\kappa}
\newcommand{\no}{\nonumber}
\newcommand{\s}{\sqrt}
\newcommand{\barg}{\bar{g}}
\newcommand{\barp}{\bar{\psi}}
\newcommand{\bara}{\bar{a}}
\newcommand{\e}{\mathbb{E}}
\newcommand{\tD}{\widetilde{\De}}
\newcommand{\bY}{\bar{Y}}

\begin{frontmatter}
\title{Self-normalized Cram\'er Type Moderate Deviations under Dependence}
\center{June 17, 2014}
\runtitle{Moderate Deviations for Studentized Statistics}

\begin{aug}
\author{\fnms{Xiaohong}
\snm{Chen}\protect\ead[label=e1]{xiaohong.chen@yale.edu}\thanksref{t1}},
\thankstext{t1}{Research supported by Cowles Foundation for Research in Economics.}
\author{\fnms{Qi-Man}
\snm{Shao}\protect\ead[label=e2]{qmshao@sta.cuhk.edu.hk}\thanksref{t2}}
\thankstext{t2}{Research partially supported by Hong Kong RGC GRF 403513 and 603710.}
and \author{\fnms{Wei Biao}
\snm{Wu}\protect\ead[label=e3]{wbwu@galton.uchicago.edu}\thanksref{t3}},
\thankstext{t3}{Research partially supported by DMS-0906073 and DMS-1106790.}

\runauthor{X. Chen, Q.M. Shao and W.B. Wu}

\affiliation{Yale University, The Chinese University of Hong Kong
and The University of Chicago}

\address{Department of Economics\\
Yale University \\
30 Hillhouse Ave., Box 208281\\
New Haven, CT 06520\\
\printead{e1}\\
\phantom{E-mail:\ } }

\address{Department of Statistics\\
The Chinese University of Hong Kong\\
Shatin, N.T.\\
Hong Kong\\
\printead{e2}\\
\phantom{E-mail:\ }}

\address{Department of Statistics\\
The University of Chicago\\
5734 S University Ave\\
Chicago, IL 60637\\
\printead{e3}\\
\phantom{E-mail:\ }}

\end{aug}

\begin{abstract}
We establish a Cram\'er-type moderate deviation result for
self-normalized sums of weakly dependent random variables, where the moment requirement is much weaker than
the non-self-normalized counterpart. The range of the moderate deviation is shown to depend on the
moment condition and the degree of dependence of the underlying
processes. We consider two types of self-normalization: the
big-block-small-block scheme and the interlacing or equal-block
scheme. Simulation study shows that the latter can have a better
finite-sample performance. Our result is applied to multiple
testing and construction of simultaneous confidence intervals for
high-dimensional time series mean vectors.
\end{abstract}



\end{frontmatter}

\section{Introduction}  \lbl{intro}
\setcounter{equation}{0}

Self-normalized sums have attracted considerable attention
recently. In comparison with their non-self-normalized
counterpart, the range of Gaussian approximation can be much wider
under same moment conditions. Let $X_i, 1 \leq i \leq n$, be
independent mean zero random variables and $S_n =\sum_{i=1}^n X_i $. Define the self-normalized
sum
\begin{equation}\label{eq:D1657}
T_n = { {S_n} \over V_n}, \mbox{ where }
 V_n^2 = \sum_{i=1}^n X_i^2.
\end{equation}
Let \begin{eqnarray*}
d_{n,\delta} = (\sum_{i=1}^n E X_i^2)^{1/2}
  / ( \sum_{i=1}^n E|X_i|^{2+\delta})^{1/(2+\delta)}.
\end{eqnarray*}
The following Cram\'er type moderate deviation result is a version of Theorems 2.1 and 2.3 of Jing, Shao and Wang (2003):

\begin{theorem} \lbl{JSW03}
Let $X_i , 1 \leq i \leq n$, be independent with $E X_i=0$, $E|X_i|^2 >0$ and
$E|X_i|^{2+\delta} < \infty$ for $0  < \delta \leq 1$ and all $i$. Then there exists an
absolute constant $A$ such that
\begin{equation}
{ P( T_n \geq x) \over 1 - \Phi(x)} = 1+ O(1) {(1+x)^{2+\delta} \over d_{n,\delta}^{2+\delta} }
\lbl{JSW03-2}
\end{equation}
for all $0 \leq x \leq  d_{n, \delta}$, and  $|O(1)| \leq A$. If in addition,
\begin{equation}
(\sum_{i=1}^n EX_i^2) \max_{1 \leq j \leq n}(E(X_j^2))^{\delta /2}
\leq \sum_{i=1}^n E|X_i|^{2+\delta}
\lbl{c-1}
\end{equation}
then
\begin{equation}
{ P(T_n \geq x) \over 1-\Phi(x)}
 = \exp \left ( O(1) {(1+x)^{2+\delta} \over d_{n,\delta}^{2+\delta} } \right )
\lbl{JSW03-1}
\end{equation}
for all $0 \leq x \leq \left( d_{n, \delta}^{2+\delta} /A\right)^{1/\delta}$, and the constant $|O(1)| \leq A$.
\newline
If $E|X_i|^2 \geq c >0$ and $E|X_i|^{2+\delta} \leq c' < \infty$ for $0  < \delta \leq 1$ and all $i$, then condition (\ref{c-1}) is automatically satisfied, and equation (\ref{JSW03-1}) holds with $d_{n,\delta} \asymp n^{\delta/(4+2\delta)}$ for all $0 \leq x \leq O(n^{1/2})$.
\end{theorem}

If in (\ref{eq:D1657}), (\ref{JSW03-1}) and (\ref{JSW03-2}) we use the
non-self-normalized version with $T_n ' = S_n / (E(V^2_n))^{1/2}$,
then the range of $x$ such that (\ref{JSW03-1}) (or (\ref{JSW03-2})) holds can be much
narrower. The moderate deviation result of type (\ref{JSW03-1}) (or (\ref{JSW03-2}))
plays an important role in statistical inference of means since in
practice one usually does not know the variance ${\rm var}(S_n) =
E (V^2_n)$. Even if the latter is known, it is still advisable to
use $T_n$, due to its wider range of Gaussian approximation. For a
comprehensive study of self-normalized sums we refer to
(\cite{DLS09}) and (\cite{SW2013}).

The main purpose of this paper is to establish a Cram\'er-type
moderate deviation result for self-normalized sums of weakly
dependent random variables under weak moment conditions. This result should be very useful for statistical inference on dependent data with fat-tailed marginal distributions. In this case due
to the dependence the self-normalized denominator $V_n$ in (\ref{eq:D1657}) is no
longer valid. In the context of resampling theory for weakly dependent
processes, block bootstrap procedures were proposed for adjusting
for dependence; see \cite{MR1707286} and \cite{MR2001447}.
However, the problem of tail Gaussian approximation is rarely
studied. In this paper we first propose two types of self-normalized sums based on the big-block-small-block scheme
and the interlacing or equal-block scheme, and then establish their associated moderate deviation theory. It is shown that, due
to the dependence, the range of Gaussian approximation is narrower
than their independent counterparts, but is still wider than their non-self-normalized ones under same polynomial moment conditions.

Although we focus on establishing a self-normalized Cram\'er-type
moderate deviation result for weakly
dependent data, our proof technique could be used to extend additional self-normalized limit theorems in Jing, Shao and Wang (2003), \cite{LS2013} and others surveyed in \cite{SW2013} from independent data to weakly dependent data with finite moments.

The rest of this paper is structured as follows. Section
\ref{dependence} introduces weakly dependent processes in terms of
$\beta$-mixing coefficients and functional dependence measures. These notions of dependence are not nested. Together they cover a large class of widely used (nonlinear) time series models.
Section \ref{main} provides two types of self-normalized sums for dependent data and derive their moderate deviation theorems. It also presents a two-sample moderate deviation extension. Section
\ref{appl} gives an application to multiple
test for high-dimensional time series mean vectors, where the test in \cite{FHY07} is generalized to the
dependence setting. Section \ref{simu} presents a simulation
study, which indicates that the self-normalized sums based on the interlacing scheme performs very well in
finite samples. All the proofs are given in Section \ref{proof}.

\section{Dependence Measures} \lbl{dependence}
\setcounter{equation}{0}

There are many different notions of temporal dependence for
general (nonlinear) time series. In this paper we focus on two
measures of dependence that have been shown to cover a large
class of time series models commonly used in statistics,
econometrics, finance and engineering.

\subsection{$\beta$-mixing}\label{sec:beta}

Many widely used nonlinear time series models can be shown to be
beta-mixing and/or strong-mixing.

Let $\{X_t\}_{t=1}^{\infty}$ be a sequence of random variables that may be non-stationary.
Let $\mathcal{I}_{-\infty }^{t}$ and $\mathcal{I}_{t+j}^{\infty }$ be $%
\sigma -$fields generated respectively by ($X_{-\infty },\cdots
,X_{t}$) and
($X_{t+j},\cdots ,X_{\infty }$). Define%
\begin{equation*}
\beta (j)\equiv \sup_{t}E\sup \{|P(B|\mathcal{I}_{-\infty
}^{t})-P(B)|:B\in \mathcal{I}_{t+j}^{\infty }\}.
\end{equation*}%
\begin{equation*}
\alpha (j)\equiv \sup_{t}\sup \{|P(A\cap B)-P(A)P(B)|:A\in \mathcal{I}%
_{-\infty }^{t},B\in \mathcal{I}_{t+j}^{\infty }\}.
\end{equation*}%
$\{X_{t}\}_{t=-\infty }^{\infty }$ is called
$\beta$-\textit{mixing} (or absolutely regular) if $\beta (j)\rightarrow 0$ as
$j\rightarrow \infty $ and is \textit{strong mixing} if $\alpha
(j)\rightarrow 0$ as $j\rightarrow \infty .$

There are alternative yet equivalent definitions of these mixing
conditions for Markov processes. For a strictly stationary Markov process $%
\left\{ X_{t}\right\}_{t=1}^{\infty }$ on a set $\Omega \subseteq \mathcal{R%
}^{d}$, let $||\phi ||_{p}^{p}=\int_{\Omega }|\phi (y)|^{p}dQ(y)$ and $%
\mathcal{T}_{t}\phi (x)=E[\phi (X_{t+1})|X_{1} =x].$
The Markov process $\left\{ X_{t}\right\} $ is said to be $\alpha -mixing$ if%
\begin{equation*}
\alpha (t)=\sup_{\phi :E[\phi (X_{t})]=0,||\phi ||_{\infty }=1}||{\mathcal{T}%
}_{t}\phi ||_{1}\rightarrow 0\text{ as }t\rightarrow \infty ;
\end{equation*}%
and the Markov process $\{X_{t}\}$ is $\beta -mixing$ if
\begin{equation*}
\beta (t)=\int \sup_{0\leq \phi \leq 1}\left\vert
\mathcal{T}_{t}\phi (x)-\int \phi dQ\right\vert dQ\rightarrow
0\text{ as }t\rightarrow \infty .
\end{equation*}%
It is well-known that $2\alpha (t)\leq \beta (t)$. See,
e.g., Bradley (2007).

The notion of $\beta -mixing$ for a Markov process is closely
related to the concept called $V-ergodicity$ (in particular
$1-ergodicity$), see e.g., Meyn and Tweedie (1993). Given a Borel
measurable function $V\geq 1$ , the Markov process $\{X_{t}\}$ is
$V-ergodic$ if
\begin{equation*}
\lim_{t\rightarrow \infty }\sup_{0\leq \phi \leq V}\left\vert \mathcal{T}%
_{t}\phi (x)-\int \phi dQ\right\vert =0\text{ , for all }x\text{;}
\end{equation*}%
the Markov process $\{X_{t}\}$ is $V-uniformly$ $ergodic$ if for
all $t\geq 0 $,
\begin{equation*}
\sup_{0\leq \phi \leq V}\left\vert \mathcal{T}_{t}\phi (x)-\int
\phi dQ\right\vert \leq cV(x)\exp (-\delta t)
\end{equation*}%
for positive constants $c$ and $\delta $. A stationary process
that is \textit{V-uniformly ergodic} will be $\beta -mixing$ with an
exponential decay rate provided that $E[V\left( X_{t}\right)
]<\infty $. This connection is valuable because one can show that
a Markov time series is beta mixing by
applying the famous drift criterion (for ergodicity): There are constants $%
\lambda \in (0,1)$ and $d\in (0,\infty )$, a norm-like function
$\Gamma
()\geq 1$ and a small set $\mathbf{K}$ such that%
\begin{equation*}
E[\Gamma (X_{t})|X_{t-1}]\leq \lambda \Gamma (X_{t-1})+d\times
1\{X_{t-1}\in \mathbf{K}\}.
\end{equation*}%
In this case, $\{X_{t}\}$ is geometric ergodic and beta mixing
with an exponential decay rate. There is also a drift criterion for
sub-geometric ergodicity or beta mixing decay at a slower than
exponential rate. See, e.g., Tong (1990), Meyn and Tweedie (1993)
and Chan and Tong (2001).

Many nonlinear time series econometrics models are shown to be
beta mixing (and hence strong mixing) via Tweedie's drift
criterion approach. See, e.g., Tong (1990) for threshold models,
Chen and Tsay (1993a, b) for functional coefficient autoregressive
models and nonlinear additive ARX models, Masry and Tj$\phi $stheim (1995) for nonlinear
ARCH, Carrasco and Chen (2002) for GARCH, stochastic
volatility and autoregressive conditional duration,
Chen, Hansen and Carrasco (2010) for diffusions, Chen, Wu and Yi
(2009) and Beare (2010) for copula-based Markov models, Douc, Moulines, Olsson and van
Handel (2011) for a large class of generalized hidden Markov
models. See Fan and Yao (2003) and Chen (2012) for additional models and references.

\subsection{Functional dependence measures}
\label{wu} Another useful dependence measure for (nonlinear) time
series is the so-called \textit{functional dependence measure};
see, e.g., Wu (2005, 2011). Suppose that $(X_t)$ is a causal
process that can be represented as
\begin{eqnarray}
X_t = G_t({\cal F}_t),
\end{eqnarray}
where $G_t(\cdot)$ is a measurable function such that $X_t$ is a
well-defined random variable, and ${\cal F}_t = \sigma (\ldots, \varepsilon_{t-1},
\varepsilon_t)$. Here $\varepsilon_{t}$, $t\in \mathbb{Z}$, are
independent and identically distributed (i.i.d.) random variables.
Let $(\varepsilon _i^\ast)_{i\in
\mathbb{Z}}$ be an iid copy of $(\varepsilon_i)_{i\in
\mathbb{Z}}$, and ${\cal F}^\ast_i = \sigma (\ldots, \varepsilon^\ast_{i-1},
\varepsilon^\ast_i)$. Hence $\varepsilon_i^\ast, \varepsilon_j, i, j\in
\mathbb{Z}$, are i.i.d. Assume that, for all $t$, $X_t$ has finite
$r$th moment, $r > 2$. Define the functional dependence measures as
\begin{equation*}
\theta_r(m) = \sup_i\| X_i - G_i(\ldots, \varepsilon_{i-m-2},
\varepsilon_{i-m-1}, \varepsilon^\ast_{i-m}, \varepsilon_{i-m+1},
\ldots,  \varepsilon_i) \|_r
\end{equation*}
and
\begin{eqnarray}
\Delta_r(m) = \sup_i \|X_i - G_i({\cal F}^\ast_{i-m},
\varepsilon_{i-m+1}, \ldots, \varepsilon_i) \|_r.
\end{eqnarray}
Note that $\theta_r(m)$ is closely related to the impulse response
function for linear processes and it can be interpreted as a
nonlinear generalization of impulse response functions. We say
that $(X_t)$ is \textit{geometric moment contraction} (GMC; see Wu
and Shao 2004) if there exists $\rho \in (0, 1)$ and $0 < \tau \le
1$ such that
\begin{eqnarray}\label{eq:gmc0}
\Delta_r(m) = O(\rho^{m^\tau})= O(e^{-a_2 m^\tau})~with~ a_2
= - \log \rho.
\end{eqnarray}
It is easily seen that (\ref{eq:gmc0}) is equivalent to $\theta_r(m) =
O(\rho_1^{m^\tau})$ for some $\rho_1 \in (0, 1)$. We emphasize
that GMC does not imply geometric $\beta$-mixing. Andrews (1984)
gave a simple AR(1) example: $X_t = (X_{t-1} + \varepsilon_t )/2$,
where $\varepsilon_t $ are i.i.d. Bernoulli(1/2). This process is
not $\alpha$-mixing (and hence not $\beta$-mixing), however it satisfies GMC (\ref{eq:gmc0}) with $\rho =1/2$ (or $a_2
= \log 2$).

\textit{Examples of GMC.} Consider the infinite order
autoregressive process
\begin{eqnarray}\label{eq:N28136}
X_{k+1} = R (\varepsilon_{k+1} ; X_k, X_{k-1}, \ldots ),
\end{eqnarray}
where $\varepsilon_k$ are i.i.d. and $R$ is a measurable function;
see Wu (2011) and Doukhan and Wintenberger (2008). Assume there
exists a nonnegative sequence $(w_j)_{j \ge 1}$ with $w_* =
\sum_{j=1}^\infty w_j < 1$ such that
\begin{eqnarray*}
\|R ( \varepsilon_0; x_{-1}, x_{-2}, \ldots )
 - R (\varepsilon_0; x'_{-1}, x'_{-2}, \ldots )\|_r
 \le \sum_{j=1}^\infty w_j |x_{-j} - x_{-j}'|.
\end{eqnarray*}
By Equations (27) and (28) in Wu (2011), since $\sum_{j=1}^\infty
w_j < 1$, (\ref{eq:N28136}) has a strictly stationary solution of
the form
\begin{eqnarray*}
X_i = G (\varepsilon_i, \varepsilon_{i-1}, \ldots),
\end{eqnarray*}
whose functional dependence measures $(\theta_k)_{k \ge 0}$
satisfies
\begin{eqnarray*}
\theta_{k+1} \le \sum_{i=1}^{k+1} w_i \theta_{k+1-i}.
\end{eqnarray*}
To obtain a bound for $\theta_k$, we define the sequence $a_k$
with $a_0 = \delta_0$, and
\begin{eqnarray}\label{eq:N28210}
a_{k+1} = \sum_{i=1}^{k+1} w_i a_{k+1-i}.
\end{eqnarray}
If $w_j$ decays sub-geometrically in the sense that, for some
$\rho \in (0, 1)$, $\tau \in (0, 1)$ and $C_0 > 0$, as $j \to
\infty$,
\begin{eqnarray}\label{eq:N28212}
w_j \sim C_0 \rho^{j^\tau}.
\end{eqnarray}
Then by elementary calculations, the recursion (\ref{eq:N28210})
has the asymptotic relation
\begin{eqnarray}\label{eq:N28214}
a_k \sim { {a_0} \over { (1-w_*)^2}} C_0 \rho^{k^\tau},
\end{eqnarray}
which entails GMC condition (\ref{eq:gmc0}). If in (\ref{eq:N28212}) $\tau = 1$,
then for some $\rho_0 \in (0, 1)$, we have
\begin{eqnarray}\label{eq:N28216}
a_k \sim { {a_0} \over { (1-w_*)^2}} C_2 \rho_0^k.
\end{eqnarray}

\section{Main results} \lbl{main}
\setcounter{equation}{0}

Let $\{X_i, i \geq 1\}$ be a sequence of random variables
satisfying
\begin{equation}
 E(X_i) =\mu =0, \ \ \ E|X_i|^{r}  \leq c_1^{r} \ \ \mbox{for all} \ \ i
 \lbl{c1}
 \end{equation}
for $r>2$ and $ c_1 < \infty$. Set $S_{k,m} =\sum_{i=k+1}^{k+m}
X_i$ and $S_n = S_{0,n}$. Assume that there exists a positive
number $c_2$ such that
\begin{equation}
 E (S_{k,m}^2) \geq c_2^2 m  \ \  \mbox{for all} \ \
  k\geq 0, \ m \geq 1. \lbl{c0}
\end{equation}
We shall assume that $\{X_i\}$ is weakly dependent which can be
either geometric $\beta$-mixing or geometric moment contracting
(GMC); see Sections \ref{mdgmixing} and \ref{mdgmc}, respectively.

For independent random variables, (\ref{eq:D1657}) is the natural
form for normalized sum. The situation is quite different when
dependence is present. There are a few ways to account for
dependence. Section \ref{bbsbns} proposes the
big-block-small-block normalized sum, while Section \ref{ins}
introduces an interlacing normalized sum. For both schemes we can
establish their moderate deviations for either geometric
$\beta$-mixing or GMC processes. Blocking technique is a common way to weaken dependence; see for example Lin and Foster (2014).

\subsection{Big-block-small-block Normalized Sum}
\lbl{bbsbns} Partition $\{ X_i, \ 1 \leq i \leq n\}$ into
consecutive big blocks and small blocks. Let $m_1 = \lfloor
n^{\alpha_1} \rfloor, m_2 = \lfloor n^{\alpha_2} \rfloor$, where
$1 > \alpha_1 \geq \alpha_2 > 0$, $k = \lfloor n/(m_1+m_2)
\rfloor$ and, for $1 \leq j \leq k$,  put
\begin{eqnarray*}
 H_{j,1} &=& \{i: \ (j-1) ( m_1+m_2)+1
      \le i \le (j-1)(m_1+m_2) + m_1\}, \cr
 H_{j,2} &=& \{i: \ (j-1) ( m_1+m_2)+m_1 +1
      \le i \le j (m_1+m_2)  \},
\end{eqnarray*}
where $H_{j,1}$ (resp. $H_{j,2}$) are large (resp. small) blocks,
and the corresponding block sums
\begin{equation}
 \left\{
 \begin{array}{ll}
 Y_{j,1} & =   \sum_{i \in H_{j,1}} X_i, \ \
   Y_{j,2} = \sum_{i \in H_{j,2}} X_i, \\
 S_{n,1} & =  \sum_{j=1}^k Y_{j,1}, \ \
   S_{n,2} =  \sum_{j=1}^k Y_{j,2}, \\
 V_{n,1}^2 & =   \sum_{j=1}^k Y_{j,1}^2, \ \
   V_{n,2}^2  =   \sum_{j=1}^k Y_{j,2}^2.
\end{array} \right.
\lbl{t12}
\end{equation}
Consider the self-normalized big-block sum
\begin{equation}
W_n = { S_{n,1} \over V_{n,1}}. \lbl{wn}
\end{equation}
Under geometric $\beta$-mixing condition \eq{beta} or GMC condition \eq{eq:gmc0},  one can easily prove
that $ W_n  \sta{d.}{\to} N(0,1) $.

If the mean is common but unknown, i.e., $E(X_i)=\mu$ for all $i$ with $\mu$ unknown, we consider the Student t-statistic
\begin{eqnarray}\label{eq:wnstar}
W^*_n = { S^*_{n,1} \over V^*_{n,1}} = { {\sum_{j=1}^k (Y_{j,1} - m_1 \mu)} \over \sqrt{\sum_{j=1}^k
(Y_{j,1} - \bar Y_1 )^2}},
\end{eqnarray}
where $\bar Y_1 = k^{-1} \sum_{j=1}^k Y_{j,1}$.

\subsection{Interlacing Normalized Sum}
\lbl{ins} A particularly interesting case for $W_n$ in (\ref{wn})
is $\alpha_1 = \alpha_2 = \alpha \in (0, 1)$. Let $m = \lfloor
n^\alpha \rfloor$, $k := \lfloor n/(2 m) \rfloor$ and
\begin{eqnarray}
H_j =\{i: \ 2m(j-1) +1 \leq i \leq 2m(j-1) + m\}, \,\,
 1 \le j \le k.
\end{eqnarray}
Note that $H_j  = H_{j, 1}$. Let $Y_j = \sum_{l \in H_j} X_l$,
$V^2 = \sum_{j=1}^k Y_j^2$ and
\begin{eqnarray}\label{eq:IIn}
I_n = { {\sum_{j=1}^k Y_j} \over V} = { {\sum_{j=1}^k Y_j} \over
\sqrt{\sum_{j=1}^k Y_j^2}},
\end{eqnarray}
which is $W_n$ in (\ref{wn}). Denote by $I_n^*$ the interlaced version of $W^*_n$ in (\ref{eq:wnstar}):
\begin{eqnarray}\label{eq:IInstar}
I^*_n = { {\sum_{j=1}^k (Y_{j} - m \mu)} \over \sqrt{\sum_{j=1}^k
(Y_{j} - \bar Y )^2}},\ \  \mbox{where} \ \ \bar Y = k^{-1} \sum_{j=1}^k Y_{j}.
\end{eqnarray}

Here we shall treat $I_n$ with $m_1
= m_2$ as a separate case since we only use half of the data from
odd blocks $\{1, 2, \ldots, m\}, \, \{1+2m, 2+ 2m, \ldots, 3m\},
\ldots$, while if $\alpha_1 > \alpha_2$ we essentially use all the
data. Our simulation study shows that the equal-block scheme has a
better finite-sample normal approximation accuracy performance,
especially when the dependence of the underlying process is
relatively strong. By excluding data $X_i$ from even blocks
$\{m+1, \ldots, 2 m\}, \, \{1+3m, \ldots, 4m\}, \ldots$, we expect
that the dependence among $(Y_1, Y_2, \ldots)$ is weaker than
block sums $(Y_{1,1}, Y_{1,2}, Y_{2,1}, Y_{2,2}, \ldots)$ as in
Section \ref{bbsbns}. The better finite-sample performance of the interlacing normalized
sum $I_n$ can be intuitively explained by the fact that, due to
the dependence, for two consecutive blocks, the second block does
not add too much new information. This is especially so when the
dependence is strong. Based on this, we shall treat it as a
separate case. A similar version can be computed if we use only
even blocks.

\subsection{Moderate Deviation under Geometric $\beta$-mixing}
\lbl{mdgmixing} Assume that there exists positive numbers  $a_1,
a_2$ and $\tau$ such that
\begin{equation}
  \beta(n) \leq a_1 e^{- a_2 n^\tau} \lbl{beta}.
\end{equation}
See Subsection \ref{sec:beta} for references of examples of
geometric $\beta$-mixing processes.

\begin{theorem} \label{th1}
Assume Conditions (\ref{c1}), (\ref{c0}) and (\ref{beta}). Let $0 < \alpha_2 \leq \alpha_1 < 1$ and $0 < \delta \leq 1, \delta <
r-2$. Then there exist finite constants $c_0, A$ depending only on $c_1/c_2, a_1, a_2, r$ and $\tau$ such that
\begin{equation}
 {{P( W_n \geq x)} \over {1 - \Phi(x)}} =
  \exp( O(1) (1+x)^{2+\delta} n^{-(1-\alpha_1)\delta/2})
  \lbl{th1a}
\end{equation}
uniformly in $ 0 \leq x \leq c_0 \min( n^{(1-\alpha_1)/2}, n^{\alpha_2 \tau /2})$, and $|O(1)|\leq A$. In
particular, we have
\begin{equation}
{ P( W_n \geq x) \over 1 - \Phi(x)} = 1+ O(1) (1+x)^{2+\delta}
n^{-(1-\alpha_1) \delta/2}, \lbl{th1b}
\end{equation}
for all $0 \le x \le c_0 \min( n^{(1-\alpha_1) \delta / (4 + 2
\delta)}, n^{\alpha_2 \tau /2})$ and $|O(1)|\leq A$.
\end{theorem}

If $\tau = 1 = \delta$ and we choose $\alpha_1= \alpha_2 = 1/2$, then \eq{th1a} yields
\begin{equation}
\ln P( W_n \geq x_n) \sim - { x_n^2/2} \lbl{th1c}
\end{equation}
as $x_n \to \infty$ and $x_n = o(n^{1/4})$. Note that when  $X_i$ are independent, Theorem \ref{JSW03}
gives a wider range of $x_n = o(n^{1/2})$.

If $\tau = 1 = \delta$ and we choose $\alpha_1 = \alpha_2= 1/4$, then \eq{th1b} implies
\begin{equation}
 { P( W_n \geq x) \over 1 - \Phi(x)} = 1+ O(1) (1+x)^{3} n^{-3/8} \to 1 \lbl{th1d}
\end{equation}
uniformly in $ 0 \leq x \leq o(n^{1/8})$. Again when  $X_i$ are independent, Theorem \ref{JSW03}
gives a wider range of $0 \le x \le o(n^{1/6})$.

\ignore{
one can choose $ \alpha_1 =
\alpha_2 = 1/4$ and the range of $x$ becomes $0 \le x \le
o(n^{1/8})$. If $X_i$ are independent, then Theorem \ref{JSW03}
gives the wider range $0 \le x \le o(n^{1/6})$.
}

\ignore{

In the big-block-small-block self-normalized sum (\ref{wn}), we
let
\begin{equation}
\alpha_1 = { \delta \over \delta + \tau ( 2+\delta) /(1+\tau
(2+\delta))}, \ \ \alpha_2 = { \alpha_1 \over 1+ \tau (2+\delta)}.
\lbl{alpha-1}
\end{equation}

\begin{theorem} \label{th1}
Assume (\ref{c1}), (\ref{c0}) and (\ref{beta}). Then we have
\begin{equation}
{ P( W_n \geq x) \over 1 - \Phi(x)} = 1+ O(1) (1+x)^{2+\delta}
n^{-(1-\alpha_1) \delta/2} \lbl{th1a}
\end{equation}
where $0 \le x \le O(n^{(1-\alpha) \delta / (4 + 2 \delta)})$.
\end{theorem}
}

In practice, it is more common to use the Student t-statistic
$\W_n$ rather than the self-normalized $W_n$. \ignore{ where
\begin{eqnarray}\label{eq:Wstar}
\W_n = { S_n  \over
 ( \sum_{j=1}^k ( Y_{j,1} - m_1 \bar{X})^2
 + ( Y_{j,2}-m_2 \bar{X})^2)^{1/2}},
\end{eqnarray}
where $\bar{X} = S_n/n$. } We have the same result for $\W_n$ in (\ref{eq:wnstar}).

\begin{theorem} \label{th2}
Let Conditions (\ref{c1}) (with unknown mean $\mu$), (\ref{c0}) and (\ref{beta}) hold. Then Results (\ref{th1a})
and \eq{th1b} also hold for $\W_n$.
\end{theorem}

\subsubsection{Moderate deviation for two-sample statistic}
The results can be easily extended to two independent sequences of
$\beta$-mixing random variables. Let $\{X_{i}^{(1)}, i \geq 1\}$
and $\{X_{i}^{(2)}, i \geq 1\}$ be two independent sequences of
$\beta$-mixing random variables with the same order of mixing decay rate and satisfy
\begin{equation}
 E(X_{i}^{(l)}) =0, \ \ \ E|X_{i}^{(l)}|^{r}  \leq c_1^{r}, \ \ l=1,2 \ \  \mbox{for all} \ \ i
 \lbl{c12}
 \end{equation}
for $r>2$ and $ c_1 < \infty$. Set $S_{k,m}^{(l)}
=\sum_{i=k+1}^{k+m} X_i^{(l)}$ and $S_n^{(l)} = S_{0,n}^{(l)}$.
Assume that there exists positive numbers  $c_2, a_1, a_2 $ and
$\tau$ such that
\begin{equation}
 E\([S_{k,m}^{(l)}]^2\) \geq c_2^2 m  \ \  \mbox{for all} \ \
  k\geq 0, \ m \geq 1, \ \ l=1,2; \lbl{c02}
  \end{equation}
 and
\begin{equation}
  \beta(n) \leq a_1 e^{- a_2 n^\tau} \lbl{beta2}
\end{equation}

Assume $n_1 \asymp n_2 \asymp n$. For $l=1,2$ we partition $\{X_{i}^{(l)}, \ 1 \leq i \leq n_l\}$ into big blocks and small
blocks. Let $m_1 = \lfloor (n_1+n_2)^{\alpha_1} \rfloor, m_2 =
\lfloor (n_1+n_2)^{\alpha_2} \rfloor$, where
$1 > \alpha_1 \geq \alpha_2 > 0$, $k_l = \lfloor n_l/(m_1+m_2)
\rfloor$ for $l=1,2$, and for $1 \leq j
\leq \max(k_1,k_2)$, put
\beq
 H_{l;j,1} &= &\{i: \ (j-1) ( m_1+m_2)+1
       \leq i \leq \min( n_l, (j-1)(m_1+m_2) + m_1)\},\\
  H_{l;j,2} & = & \{i: \ (j-1) ( m_1+m_2)+m_1 +1
       \leq i \leq \min(n_l, j (m_1+m_2))  \},
\eeq
 For $l =1, 2$,
 \beq
 Y_{j,1}^{(l)} & =  & \sum_{i \in H_{l;j,1}} X_i^{(l)}, \ \
 Y_{j,2}^{(l)} = \sum_{i \in H_{l;j,2}} X_i^{(l)}, \\
 S_{n,1}^{(l)} & = & \sum_{j=1}^{k_l} Y_{j,1}^{(l)}, \ \
 S_{n,2}^{(l)}  =  \sum_{j=1}^{k_l} Y_{j,2}^{(l)}, \\
 V_{n,1}^{(l) 2} & = &  \sum_{j=1}^{k_l} [Y_{j,1}^{(l)}]^2, \ \
 V_{n,2}^{(l)2}  = \sum_{j=1}^{k_l} [Y_{j,2}^{(l)}]^2.
 \eeq
Consider
$$
\hat{W}_n = { { 1 \over k_1} S_{n,1}^{(1)} - { 1 \over
k_2} S_{n,1}^{(2)} \over \( { 1 \over  k_1^2} V_{n,1}^{(1)2} +
{ 1 \over  k_2^2}V_{n,1}^{(2)2}\)^{1/2}}.
$$

\begin{theorem} \lbl{th2}
Under Conditions \eq{c12}, \eq{c02} and \eq{beta2}, Results \eq{th1a} and
\eq{th1b} remain valid for $\hat{W}_n$.
\end{theorem}

\ignore{

\begin{theorem} \lbl{th-m}
Let $\{X_i,  i \geq 1\}$ be m-dependent random variables with
$EX_i=0$ and assume that
$$
ES_{i,j}^2 \geq c_2 ^2 \ j
$$
Then
$$
\ln P( S_n/(V_1^2 + V_2^2)^{1/2} \geq x) \sim - x^2/2
$$
and
\end{theorem}
}

\subsection{Moderate Deviation under Geometric Moment Contraction}
\lbl{mdgmc}

In this subsection we consider time series models that satisfy the
GMC condition of Wu and Shao (2004). Similarly as (\ref{beta}),
assume that there exist $a_1, a_2 > 0$ and $0 < \tau \le 1$ such
that
\begin{eqnarray}\label{eq:gmc}
\Delta_m \le a_1 e^{-a_2 m^\tau}.
\end{eqnarray}

\begin{theorem}\label{th:mdgmc}
\begin{itemize}
\item[(1)] Assume Conditions (\ref{c1}), (\ref{c0}) and (\ref{eq:gmc}).
Let $ 0 < \alpha < 1$ and $2 < r \leq 3$. Then there exist finite constants $c_0, A$ depending only on $c_1/c_2, a_1, a_2, \alpha, r$ and $\tau$ such that $I_n$ in (\ref{eq:IIn}) satisfies the following moderate deviation theorem
\begin{equation}
 { P(I_n \ge x) \over 1- \Phi(x)}
 =  \exp\big( O(1) (1+x)^r n^{(1-r/2)(1-\alpha)}\big)  \lbl{th2a}
 \end{equation}
 for all $ 0 \leq x \leq c_0 \min( n^{(1-\alpha )/2}, n^{\alpha \tau /2})$ and $|O(1)| \leq A$. In particular, we have
 \begin{eqnarray}\label{eq:N25624}
{ {P(I_n \ge x)} \over {1-\Phi(x)}}
 = 1 + O(1) (1+x)^r n^{(1-r/2)(1-\alpha)}
\end{eqnarray}
for all $ 0 \leq x \leq c_0 \min( n^{(1-\alpha )(r-2)/(2r)}, n^{\alpha \tau /2})$ and $|O(1)| \leq A$.

 \item[(2)] If Condition (\ref{c1}) holds with unknown $\mu$ in Part (1), then results \eq{th2a} and (\ref{eq:N25624}) hold with $I^*_n$ in (\ref{eq:IInstar}).
     \end{itemize}
\end{theorem}

If we increase $\tau$ or $r$, then the range for $x$ could be wider. Let
$\tau = 1$, $r = 3$ and $\alpha = 1/4$. Then the moderate
deviation (\ref{eq:N25624}) implies (\ref{th1d}) uniformly in the range $0 \le x
\le o(n^{1/8})$. In comparison, if $\alpha_1 > \alpha_2$, then the
big-block-small-block self-normalized sum (\ref{th1b}) has a
moderate deviation with a narrower range. Note that the former
only uses half of the data. This phenomenon suggests that, when
dependence is present, one can use fewer data to achieve higher
accuracy. The latter claim is justified in a simulation study in
Section \ref{simu}.

\subsection{Small Sample Corrections}
\label{sec:ssc} In our interlacing normalized sum $I_n$ in
(\ref{eq:IIn}), if $Y_j$ are i.i.d. standard normal, then $I_n
\sim t_k$, a $t$-distribution with degrees of freedom $k$. Note
that $k \sim n^{1-\alpha} / 2$, which is much smaller than $n$. In
actual application of Theorem \ref{th:mdgmc}, instead of the
normal distribution function $\Phi$, we suggest using the $t_k$
distribution. Similar claims can be made for $I^*_n$, $W_n$ and $W^*_n$ as well. See \cite{WH2009}, \cite{DHJ11} and others for similar suggestions.

\section{Applications} \lbl{appl}
\setcounter{equation}{0}
As the result of Jing, Shao and Wang (2003) has been widely
applied in statistics and econometrics for independent data, our
results should be very useful in similar applications with
spatial dependent data and time series observations. As an
illustrative yet important application, in this section we apply
our theory to a time series extension of multiple tests of Fan,
Hall and Yao (2007).

Consider the problem of estimating the mean vector ${\bf \mu} =
(\mu_1, \ldots, \mu_p)'$ of a stationary $p$-dimensional vector
process ${\bf Z}_i = (Z_{i 1}, \ldots, Z_{i p})'$. Assume that
there exists a constant $C < \infty$ such that, for all $1 \le l
\le p$, $E (|Z_{i l}|^3) \le C$. For the dependence condition,
assume either the $\beta$-mixing condition (\ref{beta}) or the
GMC condition (\ref{eq:gmc}) holds with $\tau = 1$
uniformly for all component process $(Z_{i l})_{i \in
\mathbb{Z}}$, $l = 1, \ldots, p$. Namely there exist finite
positive constants $a_1$ and $a_2$, independent of $p$, such that
(\ref{beta}) or (\ref{eq:gmc}) holds. Assume that the dimension
$p$ satisfies
\begin{eqnarray}\label{eq:J191254}
\log p = o(n^{1/4}).
\end{eqnarray}
Let $\alpha \in (0, 1)$ be a given level. Then the upper $(\alpha
/ (2 p))$th quantile $\Phi^{-1}(1-\alpha / (2 p)) = O( (\log
p)^{1/2} ) = o(n^{1/8})$. Applying Theorem \ref{th:mdgmc}, we can
construct $1-\alpha$ simultaneous confidence intervals for
$(\mu_l)_{l = 1}^p$ via the Bonferroni procedure as
\begin{eqnarray}\label{eq:J19107}
 {{{\bar Y}_l} \over m} \pm {{\Phi^{-1}(1-\alpha / (2 p))} \over {k
 m}} \sqrt{ \sum_{j=1}^k (Y_{j l} - {\bar Y}_l)^2},
\end{eqnarray}
where ${\bar Y}_l = k^{-1} \sum_{j=1}^k Z_{j l}$, $m \asymp
n^{1/4}$, $k = \lfloor n/(2m) \rfloor$. As discussed in Subsection
\ref{sec:ssc}, the finite-sample performance can be improved if in
(\ref{eq:J19107}) we use quantiles of $t$ distributions. The
simultaneous confidence intervals in (\ref{eq:J19107}) can be used
for testing the hypothesis $H_0: \mu = \mu^\circ$, namely $\mu_1 =
\mu_1^\circ, \ldots, \mu_p = \mu_p^\circ$. We reject the null
hypothesis at level $\alpha$ if there exists one of the intervals
in (\ref{eq:J19107}) that does not include the corresponding
$\mu_l^\circ$.

\begin{example}
\label{ex:lp}{\rm (High-dimensional linear process) Let $\eta_{i
j}, i \in \mathbb{Z}, 1 \le j \le p$, be i.i.d. random variables
with mean $0$, variance $1$ and with finite $r$th moment; let
$\eta_i = (\eta_{i 1}, \ldots, \eta_{i p})'$ and
\begin{eqnarray}
{\bf Z}_i =
  (Z_{i 1}, \ldots, Z_{i p})' = \sum_{j=0}^\infty A_j \eta_{i-j},
\end{eqnarray}
where $A_j = (a_{j, k l})_{k, l \le p}$ are coefficient matrices.
Assume that there exists a constant $c < \infty$ such that, for
all $k \le p$, $\sum_{j=0}^\infty \sum_{l=1}^p a_{j, k
l}^2 \le c$. Then by Rosenthal's inequality, $\|Z_{i
k} \|^2_r \le  (r-1)  \sum_{j=0}^\infty \| \sum_{l=1}^p a_{j, k
l} \eta_{j l} \|_r^2 \le c (r-1)^2$, implying (\ref{c1}). If
additionally, as $m \to \infty$, the operator norm $\rho(A_m) = O(f^m)$ holds for some
$f < 1$, then it is easily seen that (\ref{eq:gmc}) holds, and
consequently Theorem \ref{th:mdgmc} is applicable.}
\end{example}

\begin{example}
\label{ex:nlp}{\rm (High-dimensional nonlinear process) Assume that the $p$-dimensional vector process ${\bf Z}_i$ satisfies the recursion
\begin{eqnarray}\label{eq:A261246}
{\bf Z}_i = G({\bf Z}_{i-1}, \eta_{i}),
\end{eqnarray}
where $G(\cdot, \cdot) = (G_1(\cdot, \cdot), \ldots, G_p(\cdot, \cdot))^\top$ is a measurable function and  $\eta_{i}, i \in \mathbb{Z},$ are i.i.d. random variables. We now introduce a set of sufficient conditions for the stationarity of (\ref{eq:A261246}). Let $r > 0$. Assume that there exists $L \in (0, 1)$ such that
\begin{eqnarray}\label{eq:A261247}
\sup_{ {\bf z} \not= {\bf z}'} { {\|G({\bf z}, \eta_0) - G({\bf z}', \eta_0)\|_r} \over {|{\bf z} - {\bf z}'|}} \le L
\end{eqnarray}
and $C > 0, \theta > 0$ and ${\bf z}_0$ such that
\begin{eqnarray}\label{eq:A261248}
|{\bf z}_0| + \|G({\bf z}_0, \eta_0) \|_r \le C p^\theta.
\end{eqnarray}
Following the argument in \cite{Wu2004}, we conclude that the recursion (\ref{eq:A261246}) has a stationary solution of the form ${\bf Z}_i = g ({\cal F}_i)$ satisfying $\| {\bf Z}_i \|_r \le C_1  p^\theta$ for some $C_1 > 0$, and the function dependence measure
\begin{eqnarray}\label{eq:A261255}
\|{\bf Z}_i  - {\bf Z}_i '  \|_r \le C_1 p^\theta L^i, \mbox{ where } {\bf Z}_i ' = g({\cal F}_0^*, \eta_1, \ldots, \eta_i)
\end{eqnarray}
holds for all $i \ge 0$. Then the functional dependence measure for each component series $({\bf Z}_{i j})_{i \in \mathbb{Z}}$ is also bounded by $C_1 p^\theta L^i$ for all $j = 1, \ldots, p$. A careful check of the proof of Theorem \ref{th:mdgmc} indicates that (\ref{eq:N25624}) is still valid, by noticing that if in $m$-dependence approximation, the lag $m \ge \kappa n^{1/4}$ for a sufficiently large $\kappa$, then $p^\theta L^m \le L^{m/2}$ in view of (\ref{eq:J191254}). Therefore, based on the sample ${\bf Z}_1, \ldots, {\bf Z}_n$, the intervals (\ref{eq:J19107}) are $1-\alpha$ (conservative) simultaneous confidence intervals for the mean vector $\mu = E {\bf Z}_i$.
}
\end{example}

\section{A simulation study} \lbl{simu}
\setcounter{equation}{0} In this section we shall study the
finite-sample approximation accuracy in Theorems \ref{th1} and
\ref{th:mdgmc}. Consider the AR(1) process
\begin{eqnarray}\label{eq:D11227}
X_i = \rho X_{i-1} + \varepsilon_i,
\end{eqnarray}
where $\varepsilon_i$ are i.i.d. standard normal and $|\rho| < 1$,
and the ARCH(1) process
\begin{eqnarray}\label{eq:D11236}
U_i = (a^2 + b^2 U^2_{i-1})^{1/2} \varepsilon_i,
\end{eqnarray}
where $\varepsilon_i$ are also i.i.d. standard normal and $a$ and
$b$ are real parameters with $E(b \varepsilon_0)^2 < 1$, namely
$|b| < 1$. The larger the value of $\rho$ or $b$, the stronger the
dependence of the process $(X_i)$ or $(U_i)$, respectively. We
choose 10 levels: $\rho = 0, 0.1, \ldots, 0.9$ and $b = 0, 0.1,
\ldots, 0.9$. For $\rho = 0$ (resp. $b = 0$), $(X_i)$ (resp.
$(U_i)$) are i.i.d. normal. Note that, if $b = 0.9$ then
$E(|U_i|^3) = \infty$ since $E(|b \varepsilon_0|^3) > 1$. In fact,
according to Basrak, Davis and Mikosch (2002), let $p = p(b)$ be such
that $E(|b \varepsilon_0|^p) = 1$. Then $U_i$ has finite $r$th
moment with $r \in (0, p)$, but $E(|U_i|^p) = \infty$.

Choose $n = 1000$. In $W^*_n$ of (\ref{eq:wnstar}) we let $m_1 = 43$ and
$m_2 = 7$. In the interlacing version $I_n^*$ of (\ref{eq:IInstar}) we let $m = 50$. In Theorems
\ref{th2} and \ref{th:mdgmc}, we consider 25 levels of $x$: $x =
1.6, 1.7, \ldots, 4.0$. As discussed in Section \ref{sec:ssc},
instead of the Gaussian approximation, more accuracy can be gained
if we use $t$ distribution. Assuming that the data $X_i$ are
i.i.d. standard normal. Then $I_n^*$ has $t$
distribution with degrees of freedom $n/(2m) - 1 = 9$. Besides
$I_n^*$ and $W^*_n$, we also consider the self-normalized sum
\begin{eqnarray}\label{eq:D2512}
T_n^* = { {\sum_{j=1}^{2 k} B_j} \over
 \sqrt{\sum_{j=1}^{2 k} (B_j - \bar B)^2}},
\end{eqnarray}
where $B_j = \sum_{l=1+(j-1)m}^{j m} X_l$, $\bar B = (
\sum_{j=1}^{2 k} B_j ) / (2 k)$. Note that $T_n^*$ is closely
related to the non-overlap block bootstrap scheme, and $\hat
\sigma^2 := (2 k m)^{-1} \sum_{j=1}^{2 k} (B_j - \bar B)^2$ is the
non-overlap batched mean estimate of the long-run variance
$\sigma^2 = \sum_{k \in \mathbb{Z}} {\rm cov}(X_0, X_k)$; see
Politis, Romano and Wolf (1999) and B\"uhlmann (2002).

If $X_i$ are i.i.d. standard normal, then it is easily seen that
$T_n^* \sim t_{n/(2m)-1} = t_{19}$ and also $W_n^* \sim t_{k-1} =
t_{19}$. Table 1 shows the tail probabilities based on Gaussian
and $t$ distributions. The last column presents the ratio between
Gaussian and $t_9$ tail probabilities. When $x$ becomes larger, as
expected, the Gaussian approximation becomes worse.

\begin{center}
Table 1. Tail Probabilities of Gaussian and $t$ distributions.
$\begin{array}{|ccccc|} \hline
  x & 1-\Phi &1-t_{19}& 1-t_9 & (1-t_9)/(1-\Phi)\cr
\hline
 1.6&0.05480&0.06305&0.07203&1.31446\cr
 1.7&0.04457&0.05272&0.06167&1.38389\cr
 1.8&0.03593&0.04388&0.05270&1.46660\cr
 1.9&0.02872&0.03636&0.04494&1.56509\cr
 2.0&0.02275&0.03000&0.03828&1.68247\cr
 2.1&0.01786&0.02466&0.03256&1.82257\cr
 2.2&0.01390&0.02019&0.02767&1.99017\cr
 2.3&0.01072&0.01648&0.02350&2.19130\cr
 2.4&0.00820&0.01340&0.01995&2.43353\cr
 2.5&0.00621&0.01087&0.01693&2.72654\cr
 2.6&0.00466&0.00879&0.01437&3.08271\cr
 2.7&0.00347&0.00709&0.01220&3.51801\cr
 2.8&0.00256&0.00571&0.01036&4.05315\cr
 2.9&0.00187&0.00459&0.00880&4.71520\cr
 3.0&0.00135&0.00368&0.00748&5.53981\cr
 3.1&0.00097&0.00295&0.00636&6.57421\cr
 3.2&0.00069&0.00236&0.00542&7.88146\cr
 3.3&0.00048&0.00188&0.00461&9.54639\cr
 3.4&0.00034&0.00150&0.00394&11.68395\cr
 3.5&0.00023&0.00120&0.00336&14.45115\cr
 3.6&0.00016&0.00095&0.00287&18.06411\cr
 3.7&0.00011&0.00076&0.00246&22.82270\cr
 3.8&0.00007&0.00060&0.00211&29.14637\cr
 3.9&0.00005&0.00048&0.00181&37.62668\cr
 4.0&0.00003&0.00038&0.00156&49.10493\cr
  \hline
\end{array}$
\end{center}

Tables 2, 3 and 4 show the ratios $P(T_n^* \ge x) / (1-P(t_{19}
\ge x))$, $P(I_n^* \ge x)/ (1-P(t_{9} \ge x))$ and $P(W_n^* \ge
x)/ (1-P(t_{19} \ge x))$, where the probabilities $P(T_n^* \ge
x)$, $P(I_n^* \ge x)$ and $P(W_n^* \ge x)$ are approximated by
simulating $10^6$ realizations of the AR(1) process
(\ref{eq:D11227}). As the dependence becomes stronger, namely
$\rho$ is bigger, or $x$ moves away from $0$, the moderate
deviation approximations for $T_n^*$ and $W_n^*$ become worse with
the latter being slightly better, while the interlacing normalized
sum $I_n^*$ has a relatively consistent good performance. Similar
conclusions can be made for the ARCH process (\ref{eq:D11236}).
For example, when $x=4$ and $b=0.9$, with $10^6$ repetitions, the
above three ratios are $0.65$, $0.70$ and $0.63$, respectively.
Details are omitted. In practice, we suggest using $I_n^*$.

\newpage

\begin{center} Table 2. Moderate deviation ratios $P(T_n^* \ge x) / (1-P(t_{19} \ge x))$ for the AR(1) process  (\ref{eq:D11227}) with $x=1.6, \ldots, 4.0$ and $\rho=0, \ldots, 0.9$.
$\begin{array}{|ccccccccccc|} \hline
  x & \rho=0  &0.1 & 0.2&0.3 &0.4&0.5&0.6&0.7&0.8&0.9\cr
\hline
 1.6&1.00&1.00&1.01&1.01&1.02&1.04&1.06&1.08&1.14&1.33\cr
 1.7&1.00&1.00&1.01&1.02&1.03&1.05&1.06&1.09&1.15&1.36\cr
 1.8&1.00&1.00&1.02&1.02&1.03&1.05&1.06&1.10&1.17&1.40\cr
 1.9&1.00&1.00&1.02&1.01&1.03&1.06&1.06&1.10&1.18&1.43\cr
 2.0&1.00&1.00&1.02&1.02&1.03&1.06&1.07&1.12&1.20&1.47\cr
 2.1&1.00&1.00&1.02&1.02&1.03&1.06&1.07&1.13&1.22&1.51\cr
 2.2&1.00&1.00&1.02&1.02&1.04&1.07&1.08&1.14&1.24&1.55\cr
 2.3&1.00&1.00&1.02&1.01&1.04&1.07&1.08&1.15&1.25&1.59\cr
 2.4&0.99&1.00&1.03&1.01&1.05&1.08&1.09&1.16&1.27&1.64\cr
 2.5&1.00&1.00&1.03&1.03&1.05&1.08&1.10&1.16&1.28&1.68\cr
 2.6&1.00&1.01&1.03&1.03&1.05&1.09&1.10&1.17&1.30&1.74\cr
 2.7&0.99&1.01&1.04&1.03&1.05&1.09&1.12&1.18&1.31&1.79\cr
 2.8&1.00&1.01&1.06&1.02&1.05&1.09&1.14&1.20&1.32&1.86\cr
 2.9&1.00&1.01&1.06&1.02&1.05&1.09&1.14&1.21&1.34&1.91\cr
 3.0&1.00&1.01&1.06&1.04&1.05&1.10&1.16&1.22&1.35&1.97\cr
 3.1&1.00&1.02&1.08&1.04&1.06&1.10&1.16&1.23&1.38&2.02\cr
 3.2&0.99&1.03&1.08&1.04&1.07&1.12&1.17&1.24&1.38&2.08\cr
 3.3&0.99&1.02&1.08&1.04&1.07&1.13&1.19&1.26&1.40&2.13\cr
 3.4&1.00&1.04&1.06&1.05&1.09&1.14&1.19&1.27&1.42&2.20\cr
 3.5&0.99&1.03&1.07&1.07&1.09&1.16&1.20&1.26&1.45&2.23\cr
 3.6&0.98&1.03&1.07&1.09&1.10&1.15&1.21&1.26&1.49&2.29\cr
 3.7&1.00&1.02&1.07&1.08&1.08&1.15&1.23&1.28&1.50&2.34\cr
 3.8&0.99&1.00&1.09&1.12&1.08&1.15&1.23&1.27&1.52&2.43\cr
 3.9&0.97&0.98&1.08&1.11&1.10&1.15&1.24&1.30&1.57&2.57\cr
 4.0&0.98&0.94&1.07&1.13&1.17&1.20&1.22&1.30&1.58&2.68\cr
\hline
\end{array}$
\end{center}

\newpage

\begin{center} Table 3. Moderate deviation ratios $P(I_n^* \ge x)/ (1-P(t_{9} \ge x))$ for the AR(1) process  (\ref{eq:D11227}) with $x=1.6, \ldots, 4.0$ and $\rho=0, \ldots, 0.9$.
$\begin{array}{|ccccccccccc|} \hline
  x & \rho=0  &0.1 & 0.2&0.3 &0.4&0.5&0.6&0.7&0.8&0.9\cr
\hline
 1.6&1.00&1.00&1.00&1.00&1.00&1.00&1.00&1.00&1.00&1.00\cr
 1.7&1.00&1.00&1.00&1.00&1.00&1.00&1.00&1.00&1.00&1.01\cr
 1.8&1.00&1.00&1.00&1.00&1.00&1.00&0.99&1.01&1.00&1.01\cr
 1.9&1.00&1.00&1.00&0.99&1.00&1.00&0.99&1.01&1.00&1.01\cr
 2.0&1.00&1.00&0.99&1.00&1.01&1.00&0.99&1.01&1.00&1.01\cr
 2.1&1.00&1.00&1.00&1.00&1.00&1.01&0.99&1.01&1.00&1.01\cr
 2.2&1.01&1.00&1.00&1.00&1.00&1.00&0.99&1.01&1.00&1.01\cr
 2.3&1.01&1.00&1.00&1.00&1.00&1.00&0.99&1.01&1.00&1.01\cr
 2.4&1.01&0.99&1.00&1.00&1.00&1.01&0.99&1.01&1.01&1.01\cr
 2.5&1.01&0.99&1.00&1.00&1.00&1.00&0.99&1.01&1.01&1.01\cr
 2.6&1.01&0.99&1.00&1.00&1.00&1.00&0.98&1.01&1.00&1.01\cr
 2.7&1.01&0.99&1.00&1.00&1.00&1.00&0.98&1.01&0.99&1.01\cr
 2.8&1.01&0.99&1.00&1.00&1.00&1.00&0.98&1.01&1.00&1.00\cr
 2.9&1.02&0.98&1.00&1.00&1.01&1.00&0.98&1.01&0.99&1.01\cr
 3.0&1.01&0.98&1.00&1.00&1.01&1.01&0.98&1.02&0.99&1.01\cr
 3.1&1.01&0.99&1.01&1.00&1.02&1.01&0.97&1.02&0.99&1.01\cr
 3.2&1.01&0.98&1.00&0.99&1.02&1.02&0.97&1.01&0.99&1.01\cr
 3.3&1.01&0.99&1.01&0.99&1.02&1.01&0.96&1.00&0.99&1.01\cr
 3.4&1.00&0.98&1.00&0.98&1.03&1.01&0.96&1.00&0.99&1.02\cr
 3.5&1.01&0.98&1.00&0.97&1.02&1.02&0.96&1.00&0.99&1.00\cr
 3.6&1.02&0.99&1.00&0.97&1.02&1.03&0.95&1.00&1.00&1.00\cr
 3.7&1.01&1.00&0.99&0.96&1.03&1.03&0.94&0.99&1.00&1.00\cr
 3.8&1.02&1.01&0.99&0.95&1.04&1.04&0.95&1.00&1.00&1.00\cr
 3.9&1.02&1.00&0.98&0.94&1.04&1.04&0.95&0.99&1.01&0.99\cr
 4.0&1.03&1.00&0.98&0.95&1.04&1.06&0.94&0.99&1.03&0.99\cr
\hline
\end{array}$
\end{center}

\newpage

\begin{center} Table 4. Moderate deviation ratios $P(W_n^* \ge x)/ (1-P(t_{19} \ge x))$ for the AR(1) process  (\ref{eq:D11227}) with $x=1.6, \ldots, 4.0$ and $\rho=0, \ldots, 0.9$.
$\begin{array}{|ccccccccccc|} \hline
  x & \rho=0  &0.1 & 0.2&0.3 &0.4&0.5&0.6&0.7&0.8&0.9\cr
\hline
 1.6&1.00&0.99&1.00&1.00&1.00&1.00&1.00&1.01&1.03&1.19\cr
 1.7&1.00&1.00&1.00&1.00&1.01&1.00&1.00&1.00&1.04&1.21\cr
 1.8&1.00&0.99&1.00&0.99&1.00&1.00&1.00&1.01&1.04&1.23\cr
 1.9&1.00&0.99&1.00&0.99&1.00&1.00&1.00&1.01&1.05&1.25\cr
 2.0&1.00&0.99&1.00&0.99&1.00&1.01&1.00&1.01&1.06&1.27\cr
 2.1&1.00&0.99&0.99&0.98&1.00&1.00&1.00&1.01&1.06&1.28\cr
 2.2&1.00&0.99&1.00&0.98&1.00&1.01&0.99&1.01&1.06&1.31\cr
 2.3&1.00&0.99&1.00&0.98&1.00&1.01&0.99&1.01&1.06&1.33\cr
 2.4&1.00&0.99&1.01&0.98&1.00&1.01&0.99&1.01&1.07&1.35\cr
 2.5&1.00&1.00&1.00&0.99&1.00&1.01&0.99&1.01&1.07&1.38\cr
 2.6&1.00&0.99&1.01&1.00&1.00&1.01&1.00&1.01&1.08&1.40\cr
 2.7&1.00&1.00&1.01&1.00&1.00&1.02&1.01&1.02&1.07&1.42\cr
 2.8&1.00&0.99&1.03&1.00&1.01&1.02&1.01&1.01&1.07&1.44\cr
 2.9&1.00&1.00&1.03&1.00&1.01&1.02&1.02&1.01&1.08&1.47\cr
 3.0&0.99&1.00&1.02&1.02&1.01&1.02&1.02&1.00&1.08&1.49\cr
 3.1&0.99&0.99&1.03&1.02&0.99&1.01&1.03&1.00&1.09&1.52\cr
 3.2&0.97&0.99&1.03&1.02&0.99&1.00&1.04&1.00&1.10&1.55\cr
 3.3&0.96&0.97&1.05&1.02&0.99&1.01&1.04&1.00&1.13&1.60\cr
 3.4&0.98&0.97&1.08&1.03&1.00&1.00&1.05&1.03&1.13&1.63\cr
 3.5&1.00&0.98&1.07&1.02&0.99&1.01&1.04&1.01&1.13&1.67\cr
 3.6&1.02&0.98&1.04&1.00&1.00&1.02&1.04&0.99&1.13&1.69\cr
 3.7&0.98&0.99&1.05&1.00&1.00&1.01&1.03&0.96&1.11&1.73\cr
 3.8&0.98&0.98&1.04&1.01&1.00&0.99&1.03&0.97&1.10&1.76\cr
 3.9&0.96&0.97&1.01&1.02&0.97&0.97&1.05&0.98&1.08&1.82\cr
 4.0&0.99&1.01&1.02&0.98&1.00&0.98&1.06&1.01&1.09&1.85\cr
\hline
\end{array}$
\end{center}

\section{Proofs} \lbl{proof}
\setcounter{equation}{0}

The main idea of the proof is to use $m$-dependence approximation.
For $\beta$-mixing variables, we can apply Berbee's \cite{Ber87}
theorem and convert them to independent variables. For GMC
processes, we can also use $m$-dependence approximation. Then we
apply the moderate deviation of Jing, Shao and Wang (2003) for
independent random variables.
%
%
%

Before we prove Theorem \ref{th1} we first collect some preliminary lemmas.

\begin{lemma} \lbl{berbee}

Let $\xi_i, 1 \leq i \leq n$ be a sequence of random variables on
the same probability space and define $\beta^{(i)}= \beta( \xi_i,
( \xi_{i+1}, \cdots, \xi_n))$. Then the probability space can be
extended with random variables $\txi_i$ distributed as $\xi_i$
such that $\txi_i, 1 \leq i \leq n$ are independent and
$$
P(\xi_i \not= \txi_i  \ \ \mbox{for some} \ 1 \leq i   \leq n)
\leq \beta^{(1)} + \cdots + \beta^{(n-1)}.
$$
\end{lemma}
This is Lemma 2.1 of Berbee (1987). By Theorem 4.1 in Shao and Yu (1996) we have

\begin{lemma} \lbl{l2}
Under assumptions \eq{c1} and \eq{beta},  the following holds
\begin{equation}
E|S_{k,m}|^{r'} \leq c_0 m^{r'/2} c_1^{r'}, \lbl{l2a}
\end{equation}
for any $2 \leq r' < r$, $m \geq 1$, $ k\geq 0$, where $c_0$ is a
constant depending only on $r', r, a_1, a_2$ and $\tau$.
\end{lemma}

{\bf Proof of Theorem \ref{th1}}. Clearly, \eq{l2a} and \eq{c0} yield
\begin{equation}
{ \sum_{j=1}^k E|Y_{j,1}|^{2+\delta} \over (\sum_{j=1}^k
EY_{j,1}^2 )^{(2+\delta)/2}} \leq 2 c_0 n^{-(1-\alpha_1)\delta/2}
(c_1/c_2)^{2+\delta} \lbl{bound-1}
\end{equation}

Let $\tY_j, 1 \leq j \leq k$ be
independent random variables such that $\tY_j$ and $Y_{j,1}$ have the
same distribution for each $1 \leq j \leq k$. Set
$$
\tW_n = {\sum_{j=1}^k \tY_j \over ( \sum_{j=1}^k \tY_j^2)^{1/2}}.
$$
By Lemma \ref{berbee} and $k \leq n/(2 m_2 )$, we have
\begin{equation}
|P( W_n \geq x) - P(\tW_n \geq x)| \leq k \beta(m_2) \leq a_1  \exp( - 0.5 a_2 n^{\tau \alpha_2}) . \lbl{th1-01}
\end{equation}

We next apply Theorem \ref{JSW03} to $\tW_n$. It follows from \eq{JSW03-1} that
\begin{equation}
{ P( \tW_n \geq x) \over 1 - \Phi(x)} = \exp( O(1)
(1+x)^{2+\delta} n^{-(1-\alpha_1)\delta/2}) \lbl{th1-02}
\end{equation}
for all $ 0 \leq x \leq O(1) n^{(1-\alpha_1) /2}$.

This and (\ref{th1-01}) imply that there exist finite constants $c_0, A$ depending only  on $c_1/c_2, a_1, a_2, r$ and $\tau$ such that
\begin{equation}
 {{P( W_n \geq x)} \over {1 - \Phi(x)}} =
  \exp( O(1) (1+x)^{2+\delta} n^{-(1-\alpha_1)\delta/2})
  + O(1) {{\exp(- 0.5 a_2 n^{\tau \alpha_2})} \over{1 - \Phi(x)}}
\lbl{th1a-o}
\end{equation}
uniformly in $ 0 \leq x \leq c_0 n^{(1-\alpha_1)/2}$, and $|O(1)|\leq A$. This proves \eq{th1a}.

It also follows from \eq{JSW03-2} that
\begin{equation}
{ P( \tW_n \geq x) \over 1 - \Phi(x)} = 1 + O(1)
(1+x)^{2+\delta} n^{-(1-\alpha_1)\delta/2} \lbl{th1-03}
\end{equation}
for all $ 0 \leq x \leq O(1) n^{(1-\alpha_1) \delta/(2(2+\delta))}$.
This and (\ref{th1-01}) imply \eq{th1b}.\qed

\medskip

For the proof of Theorem \ref{th:mdgmc}, we need to use the
following lemma.

\begin{lemma} \lbl{l4}
Let $\zeta_i, 1 \leq i \leq n$ be independent non-negative random
variables with $E\zeta_i^p < \infty$, where $1 < p \leq 2$. Then
for any $ 0< y < \sum_{i=1}^n E\zeta_i$
\begin{equation}
P( \sum_{i=1}^n \zeta_i \leq \sum_{i=1}^n E\zeta_i - y) \leq \exp\Big(
- { (p-1)\over 4}{  y^{p/(p-1)} \over (\sum_{i=1}^n
E\zeta_i^p)^{1/(p-1)}}\Big). \lbl{l4a}
\end{equation}
\end{lemma}

{\it Proof.} When $p=2$, \eq{l4a} is Theorem 2.19 in \cite{DLS09} with a constant $1/2$. For $1 < p \leq 2$, observing that
$$
e^{-x} \leq 1 - x + x^p \ \ \ \mbox{for} \ \ x\geq 0,
$$
we have for $t>0$
\beq
\lefteqn{
P(\sum_{i=1}^n \zeta_i \leq \sum_{i=1}^n E\zeta_i - y)}
\\
& \leq & e^{-ty + t \sum_{i=1}^n E\zeta_i} Ee^{-t \sum_{i=1}^n \zeta_i} \nn \\
& \leq & e^{-ty + t \sum_{i=1}^n E\zeta_i} \prod_{i=1}^n ( 1 - t E\zeta_i + t^p E\zeta_i^p) \nn \\
& \leq & \exp( - t y + t^p \sum_{i=1}^n E\zeta_i^p).
\eeq
Letting
$$t = \Big( { y^p \over p \sum_{i=1}^n E\zeta_i^p} \Big)^{1/(p-1)}
$$
yields
\beq
\lefteqn{
P(\sum_{i=1}^n \zeta_i \leq \sum_{i=1}^n E\zeta_i - y)}
\\
& \leq& \exp\Big( - { (p-1) y^{p/(p-1)}  \over p^{p/(p-1)} (\sum_{i=1}^n E\zeta_i^p)^{1/(p-1)} } \Big) \\
& \leq & \exp\Big( - { (p-1) y^{p/(p-1)}  \over 4 (\sum_{i=1}^n E\zeta_i^p)^{1/(p-1)} } \Big),
\eeq
as desired.\qed

\medskip

{\bf Proof of Theorem \ref{th:mdgmc}.} Recall (\ref{eq:IIn}) for
$Y_j$, $1 \le j \le k$. Let
\begin{eqnarray*}
\tilde Y_j = E(Y_j | \varepsilon_l, 2mj-3m+1 \le l \le 2mj - m),
\end{eqnarray*}
and
\begin{eqnarray*}
\tilde I_n = { {\sum_{j=1}^k \tilde Y_j} \over \tilde V}, \mbox{
where }
 \tilde V^2 = \sum_{j=1}^k \tilde Y^2_j.
\end{eqnarray*}
Note that $\tilde Y_j$ are independent, and by (\ref{eq:gmc}),
\begin{eqnarray}\label{eq:D4140}
\|Y_j - \tilde Y_j \|_r \le m a_1 e^{-a_2 m^\tau}.
\end{eqnarray}
Under (\ref{eq:gmc}), since $X_l  = \sum_{i=0}^\infty {\cal P}_{l-i} X_l$, where ${\cal P}_k \cdot = E(\cdot | {\cal F}_k) -  E(\cdot | {\cal F}_{k-1})$, we have by Burkholder's (1988) martingale inequality that
\begin{eqnarray*}
\| Y_j \|_r &=& \|  \sum_{i=0}^\infty \sum_{l \in H_j} {\cal P}_{l-i} X_l \|_r \cr
 &\le& \sum_{i=0}^\infty \| \sum_{l \in H_j} {\cal P}_{l-i} X_l \|_r \cr
 &\le& \sum_{i=0}^\infty (r-1)^{1/2}  (\sum_{l \in H_j}  \| {\cal P}_{l-i} X_l \|^2_r)^{1/2} \cr 
 &\le& (r-1)^{1/2}  \sum_{i=0}^\infty  (m \theta^2_r(i) )^{1/2} 
 =  c_3 m^{1/2}, 
\end{eqnarray*}
where $c_3 =  (r-1)^{1/2}  \sum_{i=0}^\infty  \theta_r(i) < \infty$. By condition (\ref{c0}) and
(\ref{eq:D4140}), there exists a constant $c_5 > 0$ such that $E
\tilde V^2 \ge c_5 n$. By Lemma \ref{l4} with $p = r/2$, $\zeta_j
= \tilde Y^2_j$ and $y = c_5 n / 2$, we have by elementary
calculations that
\begin{eqnarray}\label{eq:D4135}
P(\tilde V^2 \ge c_5 n / 2) \ge 1 - \exp(- c_6 k)\geq 1- \exp( -c'_6 n^{1-\alpha})
\end{eqnarray}
for some constants $c_6 , c'_6 > 0$. Also (\ref{eq:D4140}) and $m \asymp n^{\alpha}$ imply
\begin{eqnarray*}
 P(|Y_j - \tilde Y_j| \ge n^{-9})
 \le n^{9 r} m^r a^r_1 e^{-r a_2 m^\tau} = O(1) \exp(-r a_2 n^{\tau \alpha} /2).
\end{eqnarray*}
Hence there exist $c_7, c_8 > 0$ such that
\begin{eqnarray}\label{eq:D4153}
P(|I_n - \tilde I_n| \ge n^{-2}, \tilde V^2 \ge c_5 n)
 &\le& c_7 n^{c_8} e^{-r a_2 m^\tau} \cr
 &=& O(1) \exp(-r a_2 n^{\tau \alpha} /2).
\end{eqnarray}
Observe that
\begin{eqnarray}\label{eq:D4156}
\max_{0 \le x \le n} \left| {{1-\Phi(x)} \over{1-\Phi(x \pm
n^{-2})}}-1 \right| = O(n^{-1}).
\end{eqnarray}
For $0 \leq x \leq c_0 n^{\min((1-\alpha), \tau \alpha)/2}$ with a small constant $c_0>0$, it is easy  to see that
$$
\exp( -c'_6 n^{1-\alpha}) + \exp(-r a_2 n^{\tau \alpha} /2) = o(1) (1-\Phi(x)) \exp\Big( O(1) {{(1+ x)^r} \over { n^{(1-\alpha) (r-2)/2}}} \Big).
$$
\ignore{
Let $x_n = n^{\tau \alpha / 2}$. Then $x_n^2 \sim m^\tau$.
Therefore we can choose $c_0 > 0$ such that $n^{c_8} e^{-r a_2
m^\tau} = o(1-\Phi(c_0 x_n))$. Since $\tau \alpha < 1 - \alpha$,
we have $x_n^2 = o(k)$.
}

Applying Theorem \ref{JSW03} to $\tilde I_n$, we have,
for some constant $c_4 > 0$, that
\begin{eqnarray}\label{eq:D4125}
{ {P(\tilde I_n \ge x)} \over{1-\Phi(x)}} = \exp\Big( O(1){{(1+x)^r} \over {
k^{r/2-1}}} \Big )
\end{eqnarray}
for $0 \le x \le c_4 k^{1/2}$.
Hence \eq{th2a} follows from
(\ref{eq:D4125}), (\ref{eq:D4135}), (\ref{eq:D4153}) and
(\ref{eq:D4156}) with elementary calculations. (\ref{eq:N25624}) follows similarly.
\qed

\ignore{

\begin{lemma} \lbl{l4}
Let $\zeta_i, 1 \leq i \leq n$ be independent non-negative random
variables with $E\zeta_i^p < \infty$, where $1 < p \leq 2$. Then
for any $ 0< y < \sum_{i=1}^n E\zeta_i$
\begin{equation}
P( \sum_{i=1}^n \zeta_i \leq \sum_{i=1}^n E\zeta_i - y) \leq \exp(
- { (p-1)\over 4}{  y^{p/(p-1)} \over (\sum_{i=1}^n
E\zeta_i^p)^{1/(p-1)}}). \lbl{l4a}
\end{equation}
\end{lemma}

The inequality should be known. The proof is based on the
observation that
$$
e^{-t} \leq 1- t + t^{p}  \ \ \mbox{for} \ \ t \geq 0
$$

\medskip

\begin{lemma} \lbl{th1-l1}
\begin{equation}
P(V_1^2 \leq c_2^2 n/2) \leq C \exp(-a_2 n^{\tau \alpha_2}) +
\exp( - C  n^{1-\alpha_1} (c_2/c_1)^{(4+2\delta)/\delta})
\lbl{th1-l1a}
\end{equation}
\end{lemma}

\proof  By Lemma \ref{berbee}, we have \beq P(V_1^2 \leq c_2^2
n/4)
&\leq &   k \beta(m_2) + P( \sum_{j=1}^k \tilde{Y}_j^2 \leq  c_2^2 n/4)\\
& \leq & C \exp(-a_2 n^{\tau \alpha_2}) + P( \sum_{j=1}^k
\tilde{Y}_j^2 \leq c_2^2 n/4). \eeq

\blue{need to assume $x^2 < a_2 n^{\tau \alpha_2/2}$}

Noting that by \eq{c0} and \eq{l2a},
$$
\sum_{j=1}^k  E\tilde{Y}_{j,1}^2  \geq (k-1) c_2^2 m_1 \geq n
c_2^2 /2,
$$
$$
\sum_{j=1}^k E|\tilde{Y}_{j,1}|^{2+\delta} \leq C \,
c_1^{2+\delta} \, n^{1+\alpha_1 \delta/2}
$$
So by Lemmas \ref{l4} and \ref{l2}, we get \beq \lefteqn{ P(
\sum_{j=1}^k \tilde{Y}_{j,1}^2 \leq a_1 c_2^2 n)}
\\
& \leq & \exp( - C  n^{1-\alpha_1} (c_2/c_1)^{(4+2\delta)/\delta})
\eeq Thus,

\blue{ need to assume that $ x^2 < C n^{1-\alpha_1}
(c_2/c_1)^{(4+2\delta)/\delta}/2 $}

\medskip

\subsection{Proof  of Theorem \ref{th1}}

We first show that
\begin{equation}
{ P( T_1 +T_2  \geq x (V_1^2 + V_2^2)^{1/2}) \over 1- \Phi(x)}
\leq 1+ O(1) \lbl{th1-0}
\end{equation}
\red{ Let $(1+x)^{1+\delta} n^{-(1-\alpha_1)\delta/2} \leq y \leq
(x+1)/2$}. Without  loss of generality, assume $c_2=1$.

Notice that \beqn \lefteqn{
 P( T_1 +T_2  \geq x (V_1^2 + V_2^2)^{1/2}) } \nn \\
 & \leq & P( T_1 \geq (x-y) (V_1^2 + V_2^2)^{1/2})
 + P( T_2 \geq y (V_1^2 + V_2^2)^{1/2}, T_1 +T_2  \geq x (V_1^2 + V_2^2)^{1/2}) \nn \\
 & \leq & P( T_1 \geq (x-y) V_1) + P(V_1^2 \leq c_2^2 n /4)+H_1, \lbl{th1-1}
 \eeqn
 where
 $$
 H_1 =  P(T_2 \geq y c_2  n^{1/2}/2, \ T_1 +T_2  \geq x (V_1^2 + V_2^2)^{1/2}).
 $$
By Lemma \ref{berbee}, \eq{JSW03-1} and \eq{bound-1}, we have
\beqn \lefteqn{
P(T_1\geq (x-y) V_1)} \nn \\
 & \leq  & k \beta(m_2) + (1-\Phi(x-y)) \exp(O(1) (1+x)^{2+\delta} n^{-(1-\alpha_1) \delta} (c_2/c_1)^{2+\delta}) \nn \\
 & \leq & C \exp(- a_2 n^{\tau \alpha_2}) + (1-\Phi(x-y)) \exp(O(1) (1+x)^{2+\delta} n^{-(1-\alpha_1) \delta} (c_2/c_1)^{2+\delta})
 \lbl{th-1-1-1}
 \eeqn
 for $ x^{\delta} \leq n^{(1-\alpha_1)\delta} (c_1/c_2)^{2+\delta}/C$.
 It follows from Lemma \ref{th1-l1} that
 \begin{equation}
P(V_1^2 \leq c_2^2 n/4) \leq C \exp(-a_2 n^{\tau \alpha_2}) +
\exp( - C  n^{1-\alpha_1} (c_2/c_1)^{(4+2\delta)/\delta})
\lbl{th1-l1b}
\end{equation}

Therefore, by Lemma \ref{th1-l2} below, we obtain \beqn \lefteqn{
 P( T_1 +T_2  \geq x (V_1^2 + V_2^2)^{1/2}) } \nn \\
& \leq & C \exp(- a_2 n^{\tau \alpha_2}) + (1-\Phi(x-y)) \exp(O(1) (1+x)^{2+\delta} n^{-(1-\alpha_1) \delta} (c_2/c_1)^{2+\delta})\nn \\
& & + \exp( - C  n^{1-\alpha_1} (c_2/c_1)^{(4+2\delta)/\delta})
\lbl{th1-3} \eeqn

\begin{lemma}\lbl{th1-l2}
\begin{equation}
H_1 \leq  \lbl{th1-l2a}
\end{equation}
\end{lemma}

\proof To estimate $H_2$, we truncate ${Y}_{j,2}$. For the bounded
variables, we can apply the Bennett-Hoeffding exponential
inequality, and for the unbounded term we use the following
inequality (see (5.7) in Jing, Shao and Wang (2003)): for $c\geq
0, x > 1$
\begin{equation}
\{ s+ t \geq x \sqrt{ c + t^2}\} \subset \{ s \geq (x^2 -1)^{1/2}
\sqrt{c}\} \lbl{01}
\end{equation}
to derive a recursive formula. Let $m_3 = n^{\alpha_3} /(1+x)$,
 \beq
 Y_{j,2,1} & =  & {Y}_{j,2} I\{  |{Y}_{j,2}| \leq m_3 \}, \ \
  Y_{j,2,2} = {Y}_{j,2} I\{  |{Y}_{j,2}| >  m_3\}, \\
 T_{2,1} & =& \sum_{j=1}^k Y_{j,2,1}, \ \ T_{2,2}  = \sum_{j=1}^k Y_{j,2,2}.
 \eeq
 Then
 \beqn
 H_1 & \leq & P( T_{2,1} \geq y V_1) +
 \sum_{j=1}^k P( T_1 +T_2  \geq x (V_1^2 + V_2^2)^{1/2}, \ |{Y}_{j,2}| > m_3)
\lbl{th1-4}
 \eeqn

 \begin{equation}
 \sum_{i=1}^k  EY_{j,2,1}^2 \leq c_1^2 k m_2 \leq c_1^2 n^{1+\alpha_2 -\alpha_1} \lbl{th1-l2-1}
 \end{equation}

 By \eq{th1-l1b}, Lemma \ref{berbee} and the  Bennett-Hoeffding exponential inequality

 \beqn
 P(T_{2,1} \geq y V_1) & \leq & P(V_1^2 \leq n/4) +
 P(T_{2,1} \geq y n^{1/2}/2) \nn \\
 & \leq & C \exp(-a_2 n^{\tau \alpha_2}) + \exp( - C  n^{1-\alpha_1} (c_2/c_1)^{(4+2\delta)/\delta})\nn \\
 & & + k \beta(m_1)  + \exp( - { y^2 n  \over 8 c_1^2 n^{1-\alpha_1 + \alpha_2}})
 + \exp( - { y n^{1/2} \over m_3}) \nn \\
 & \leq &
  C \exp(-a_2 n^{\tau \alpha_2}) + \exp( - C  n^{1-\alpha_1} (c_2/c_1)^{(4+2\delta)/\delta})\nn \\
 & &  + \exp( - C y^2 n^{\alpha_1 - \alpha_2})
 + \exp( - { y n^{1/2} \over m_3}) \lbl{th1-l2-2}
 \eeqn
 To estimate the second term on the right hand side of \eq{th1-4},
 let $T_1^{(j)} = \sum_{i \not = j, j+1} Y_{i,1}$ be the sum without containing the neighborhood of
 $Y_{j,2}$. Similarly, let $T_2^{(j)} = \sum_{i \not = j} Y_{i,2}$ and define
  $V_1^{(j)}$ and $V_2^{(j)}$ accordingly. Then, first applying \eq{01} and then Lemma \ref{berbee}
(considering $ \{Y_{i,1}, Y_{i,2}, i \leq j-1\}, Y_{j,2}, \{
Y_{i,1}, i \geq j+2, Y_{i,2}, i \geq j+1\}$) yields \beqn
\lefteqn{
P( T_1 +T_2  \geq x (V_1^2 + V_2^2)^{1/2}, |{Y}_{j,2}| > m_3) } \nn \\
& \leq & P( T_1^{(j)} +T_2^{(j)}  \geq (x^2-3)^{1/2} (V_1^{(j)2} +
V_2^{(j)2})^{1/2},
|{Y}_{j,2}| > m_3)  \nn \\
& \leq & 2 \beta(m_1) + P(|{Y}_{j,2}| > m_3) H_{1,j} \lbl{th1-5}
\eeqn where
\begin{equation}
H_{1,j} = P( T_1^{(j)*} +T_2^{(j)*}  \geq (x^2-3)^{1/2}
(V_1^{(j)*2} + V_2^{(j)*2})^{1/2}) \lbl{th1-6}
\end{equation}
and
$$T_1^{(j)*} + T_2^{(j)*} =
(\sum_{i \leq j-1} ( Y_{i,1}+Y_{i,2}))^{\tilde{}} + (\sum_{i >
j+1} Y_{i,1} + \sum_{i\geq j+1} Y_{i,2})^{\tilde{}}.
$$
Observe that \beqn P(|{Y}_{j,2}| > m_3) \leq
E|{Y}_{j,2}|^{2+\delta} /m_3^{2+\delta} \leq C n^{\alpha_2
(2+\delta)/2} c_1^{2+\delta} / m_3^{2+\delta} \lbl{th1-7} \eeqn

\subsection{Proof of Theorem \ref{th2}}

Observe that
$$
\sum_{j=1}^m ( Y_{j,1} -n_1 \bar{X})^2 + ( Y_{j,2}-n_2 \bar{X})^2
= \sum_{j=1}^m (Y_{j,1}^2 + Y_{j,2}^2) -
$$

}

\end{document}